\documentclass[10pt,draft,twoside,leqno,a4paper]{article}
\usepackage{amsmath,amsthm,amsfonts,amssymb,verbatim}
\newtheorem{Theorem}{Theorem}
\newtheorem{Lemma}{Lemma}
\newtheorem{Corollary}{Corollary}

\newtheorem{Example}{Example}
\newcommand{\R}{\mathbb R}

\newcommand{\de}{\partial}
\renewcommand{\d}{{\mathrm d}}

\newcommand{\dist}{\mathrm{dist}}

\newcommand{\bra}{\langle}
\newcommand{\ket}{\rangle}

\newcommand{\onabla}{\overline\nabla}
\newcommand{\balpha}{\bar\alpha}

\newcommand{\xo}{x_0}
\newcommand{\ro}{r_0}

\newcommand{\urho}{u_\rho}
\newcommand{\utheta}{u_\theta}

\newcommand{\utr}{\frac{u_\theta}{\rho}}

\newcommand{\Sr}{{S_r}}
\newcommand{\aij}{a_{ij}}

\newcommand{\wmin}{\underline w}
\newcommand{\es}{\mathrm{ess\,sup}}
\newcommand{\gxo}{g_{x_0}}
\newcommand{\thetao}{\theta_0}
\renewcommand{\d}{\,{\mathrm d}}
\renewcommand{\onabla}{\overline\nabla}
\begin{document}
\title{A sharp H\"older estimate for elliptic equations in two variables}
\author{Tonia Ricciardi\thanks{
Partially supported by 
Regione Campania L.R.~5/02 and by 
the MIUR National Project {\em Variational Methods and
Nonlinear Differential Equations}}\\
{\small Dipartimento di Matematica e Applicazioni}\\
{\small Universit\`a di Napoli Federico II}\\
{\small Via Cintia, 80126 Naples, Italy}\\
{\small fax: +39 081 675665}\\
{\small\tt{ tonia.ricciardi@unina.it}}\\
}
\date{}
\maketitle
\begin{abstract}
We prove a sharp H\"older estimate for solutions of linear 
two-dimension\-al, divergence form elliptic equations with measurable coefficients,
such that the matrix of the coefficients is symmetric and has {\em unit determinant}.
Our result extends some previous work by Piccinini and Spagnolo~\cite{PS}. 
The proof relies on a 
sharp Wirtinger type inequality.
\end{abstract}
\begin{description}
\item {\textsc{Key Words:}} linear elliptic equation, measurable coefficients,
H\"older regularity
\item {\textsc{MSC 2000 Subject Classification:}} 35J60
\end{description}
\section{Introduction}
\label{sec:intro}
Let $\Omega$ be a bounded open subset of $\R^n$ and let
$u\in H_{\mathrm{loc}}^1(\Omega)$ be a weak solution to the linear elliptic equation 
in divergence form
\begin{equation}
\label{elliptic}
(\aij u_{x_i})_{x_j}=0\qquad\qquad\mathrm{in\ }\Omega,
\end{equation}
where $\aij $, $i,j=1,2,\ldots,n$, are bounded measurable functions in $\Omega$
satisfying the ellipticity condition
\begin{align}
\label{ellipticbound}
&\lambda|\xi|^2\le \aij (x)\xi_i\xi_j\le\Lambda|\xi|^2,
\end{align}
for all $x\in\Omega$, for all $\xi\in\R^n$ and for some $0<\lambda\le\Lambda$.
It is well known (see, e.g., \cite{DG,Gi,GT,Na} and references therein)  
that solutions to \eqref{elliptic} are 
$\alpha$-H\"older continuous in $\Omega$ for some $0<\alpha<1$. More precisely,
for every compact subset $K\subset\subset\Omega$ there holds 
\begin{equation}
\label{holder}
\sup_{x,y\in K, x\neq y}\frac{|u(x)-u(y)|}{|x-y|^\alpha}<+\infty.
\end{equation}
Moreover, estimates for the H\"older exponent $\alpha$ may be obtained,
depending on the ellipticity constant $L=\Lambda/\lambda$ only.
In this note we consider the case $n=2$, which unless otherwise stated
we assume henceforth.
For this case, sharp H\"older estimates were obtained by
Piccinini and Spagnolo in \cite{PS} under the symmetry assumption 
\begin{equation}
\label{symmetric}
\aij (x)=a_{ji}(x)\qquad\qquad\mathrm{in\ }\Omega.
\end{equation}
We collect in the following theorem some results from \cite{PS} which are relevant to our
considerations.
\begin{Theorem}[Piccinini and Spagnolo~\cite{PS}]
\label{thm:ps}
Let $u\in H_{\mathrm{loc}}^1(\Omega)$ be a weak solution to \eqref{elliptic}.
\begin{itemize}
\item[(i)]
If the coefficients $\aij $
satisfy \eqref{ellipticbound} and \eqref{symmetric}, then $\alpha=L^{-1/2}$. 
\item[(ii)]
If the coefficients $\aij $
satisfy \eqref{ellipticbound}, \eqref{symmetric} and if furthermore 
$\aij =a(x)\delta_{ij}$ for some
measurable function $a=a(x)$ satisfying 
$\lambda\le a(x)\le\Lambda$ for all $x\in\Omega$, then
$\alpha=\frac{4}{\pi}\arctan L^{-1/2}$.
\end{itemize}
These values of $\alpha$ are sharp.
\end{Theorem}
In particular, Theorem~\ref{thm:ps}--(ii) implies that the optimal value of $\alpha$
increases by restricting the matrices $A$ to the class of isotropic matrices
$A=a(x)I$. Thus, it is natural to seek other classes of
matrices $A$ for which the corresponding optimal value of $\alpha$ may be improved.
Our main result in this note is to show that if $A=(\aij)$ satisfies 
\eqref{ellipticbound}, \eqref{symmetric} and if furthermore $A$ has {\em unit determinant}, 
namely, if
\begin{equation}
\label{detA}
\det A(x)=a_{11}(x)a_{22}(x)-a_{12}(x)a_{21}(x)\equiv1
\qquad\qquad\mathrm{in\ }\Omega,
\end{equation}
then a more
accurate {\em integral} characterization of $\alpha$ may be obtained.
Moreover, our result is {\em sharp} within the class of matrices satisfying 
\eqref{ellipticbound}, \eqref{symmetric} and \eqref{detA}.
It should be mentioned that condition \eqref{detA} 
is relevant in the context of quasiharmonic fields, 
see \cite{IS}. 
\begin{Theorem}[Main result]
\label{thm:main}
Let $A=(\aij)$ satisfy 
\eqref{ellipticbound}, \eqref{symmetric} and \eqref{detA} in 
$\Omega$ and let $u$ satisfy \eqref{elliptic}.
Then $u$ is $\alpha$-H\"older continuous in $\Omega$,
with $\alpha$ given by  
\begin{equation}
\label{alpha}
\alpha=2\pi\left(\sup_{\xo\in\Omega}\es_{0<r<\dist(\xo,\de\Omega)}\int_{|\xi|=1}
\aij (\xo+r\xi)\xi_i\xi_j\d\sigma_\xi\right)^{-1}.
\end{equation}
\end{Theorem}
We note that under assumption \eqref{detA}, we may choose $\lambda=1/\Lambda$
in \eqref{ellipticbound}
and therefore the ellipticity constant takes the value
$L=\Lambda^2$. Hence, the Piccinini-Spagnolo estimate
in Theorem~\ref{thm:ps}--(i) yields in this case $\alpha=\Lambda^{-1}$.
On the other hand, recalling that
$\Lambda=\sup_{x\in\Omega}\sup_{|\xi|=1}\aij (x)\xi_i\xi_j$,
it is clear that 
$\alpha\ge\Lambda^{-1}$.
\par
Theorem~\ref{thm:main} implies the following
\begin{Corollary}
\label{cor:main}
Let $A=(a_{ij})$ satisfy
\eqref{ellipticbound}, \eqref{symmetric} and \eqref{detA} 
and let $u$ satisfy \eqref{elliptic}.
Then the least upper bound for the admissible values of the H\"older exponent 
of $u$ is given by
\begin{equation}
\label{balpha}
\bar\alpha=2\pi\left(\sup_{\xo\in\Omega}
\inf_{0<\ro<\dist(\xo,\de\Omega)}\es_{0<r<\ro}\int_{|\xi|=1}
\aij (\xo+r\xi)\xi_i\xi_j\d\sigma\right)^{-1}.
\end{equation}
\end{Corollary}
\par
Theorem~\ref{thm:main} is \textit{sharp}, in the sense of the following 
\begin{Example}[Sharpness]
\label{ex:sharp}
Let $\Omega=B$ the unit ball in $\R^2$, let $\theta=\arg x$ and let 
\begin{equation}
\label{Adef}
A(x)=
\frac{1}{k(\theta)}I+\left(k(\theta)-\frac{1}{k(\theta)}\right)
\frac{x\otimes x}{|x|^2}
\qquad\mathrm{in}\ B\setminus\{0\}.
\end{equation}
where 
$k:\R\to\R^+$ is a $2\pi$-periodic measurable function
bounded from above and away from $0$. 
Then $\det A(x)\equiv1$. 
By a suitable choice of $k$, we may obtain that
\begin{equation}
\label{exbalpha}
\balpha=2\pi\left(\int_0^{2\pi}k\right)^{-1}.
\end{equation}
On the other hand the function $u\in H^1(B)$ defined by
\begin{equation}
\label{uex}
u(x)=|x|^{\balpha}\cos\left(\balpha\int_0^{\arg x} k\right)
\end{equation}
satisfies equation \eqref{elliptic} with $A$ given by \eqref{Adef}. 
Clearly, its H\"older exponent is exactly $\balpha$.
\end{Example}
A verification
of Example~\ref{ex:sharp} is provided in the Appendix.
In a forthcoming note, we will show that the functions defined by \eqref{uex}
are also of interest in the context of quasiconformal maps.
The remaining part of this note is devoted to the proof
of Theorem~\ref{thm:main}.
\subsubsection*{Notation}
For a fixed $\xo\in\Omega$,
let $x=\xo+\rho e^{i\theta}$, be the polar coordinate transformation
centered at $\xo$.
We denote $u_{x_i}=\de u/\de x_i$,
$\urho=\de u/\de\rho$, $u_\theta=\de u/\de\theta$.
Following the notation in \cite{PS}, we denote $\nabla u=(u_{x_1},u_{x_2})$ and
$\onabla u=(\urho,\utr)$. 
We denote by $J(\theta)$ the rotation matrix
\[
J(\theta)=\left(
\begin{matrix}
&\cos\theta&&-\sin\theta\\
&\sin\theta&&\quad\cos\theta
\end{matrix}
\right).
\]
In this notation, we have 
\begin{equation}
\label{Jnabla}
\nabla u=J(\theta)\onabla u.
\end{equation}
Finally, we denote by $\bra\cdot,\cdot\ket$ the Euclidean scalar product
in $\R^2$.
\section{Proof of Theorem~\ref{thm:main}}
\label{sec:proof}
In the spirit of \cite{PS}, a key ingredient in the proof of Theorem~\ref{thm:main}
is a sharp
Wirtinger type inequality.
\begin{Lemma}[Sharp Wirtinger type inequality]
Let $a:\R\to\R$ be a $2\pi$-periodic measurable function which is bounded from above
and away from 0.
Then,
\begin{equation}
\label{wirtinger}
\int_0^{2\pi}aw^2\le
\left(\frac{1}{2\pi}\int_0^{2\pi}a\right)^2
\int_0^{2\pi}\frac{1}{a}(w')^2
\end{equation}
for every $w\in H^1_{\mathrm{loc}}(\R)$ such that $w$ is $2\pi$-periodic
and satisfies $\int_0^{2\pi}aw=0$.
Equality in \eqref{wirtinger} is attained if and only if $w$ is of the form
\begin{equation}
\label{minimizer}
w(\theta)=C\cos\left(\frac{2\pi}{\int_0^{2\pi}a}\int_0^\theta a
+\varphi\right)
\end{equation}
for some $C\neq0$ and $\varphi\in\R$.
\end{Lemma}
\begin{proof}
Let 
\[
X=\left\{w\in H^1_{\mathrm{loc}}(\R):\ 
w\ \mathrm{is}\ 2\pi\mathrm{-periodic}\right\}
\]
and consider the functional $I$ defined by \[
I(w)=\frac{\int_0^{2\pi}\frac{1}{a}(w')^2}{\int_0^{2\pi}aw^2}
\]
for all
$w\in X\setminus\{0\}$.
We equivalently have to show that 
\[
\inf\left\{I(w):\ w\in X\setminus\{0\}, \int_0^{2\pi}aw=0\right\}
=\left(\frac{2\pi}{\int_0^{2\pi}a}\right)^2
\]
and that the infimum is attained exactly on functions of the form
\eqref{minimizer}.
By Sobolev embeddings and standard compactness arguments, there
exists $\wmin\in X\setminus\{0\}$ satisfying $\int_0^{2\pi}a\wmin=0$ and such that 
\[
I(\wmin)=\inf\left\{I(w):\ w\in X, \int_0^{2\pi}aw=0\right\}>0.
\]
By the minimum property of $\wmin$ we have that 
$I'(\wmin)\left(\phi-(2\pi)^{-1}\int_0^{2\pi}a\phi\right)=0$
for all $\phi\in X$. By properties of $\wmin$ it follows that
\[
\int_0^{2\pi}\frac{1}{a}\wmin'\phi'-I(\wmin)\int_0^{2\pi}a\wmin\phi=0
\qquad\forall\phi\in X.
\]
Let $\Theta=\int_0^\theta a$, $w(\theta)=W(\Theta)$, $\phi(\theta)=\Phi(\Theta)$.
Then $W,\Phi\in H^1_{\mathrm{loc}}(\R)$ 
are $\left(\int_0^{2\pi}a\right)$-periodic and
\[
\int_0^{\int_0^{2\pi}a}W'(\Theta)\Phi'(\Theta)\d\Theta
=I(\wmin)\int_0^{\int_0^{2\pi}a}W(\Theta)\Phi(\Theta)\d\Theta.
\]
Equivalently, $W$ is a weak solution for
\[
-(W')'=I(\wmin)W.
\]
It follows that $I(\wmin)=(2\pi/\int_0^{2\pi}a)^2$ and
\[
W(\Theta)=C\cos(\frac{2\pi}{\int_0^{2\pi}a}\Theta+\varphi)
\]
for some $C\neq0$ and $\varphi\in\R$. Recalling the definition of $W$
and $\Theta$, we conclude the proof.
\end{proof}
We can now provide the
\begin{proof}[Proof of Theorem~\ref{thm:main}]
For every $0<r<\dist(\xo,\de\Omega)$ let
\[
\gxo(r)=\int_{|x-\xo|<r}\aij u_{x_i}u_{x_j}\d x.
\]
It is well-known (see, e.g., \cite{Gi,GT,IM}) that 
a sufficient condition for \eqref{holder} 
is that the function $G_{\xo}(r)=r^{-2\alpha}\gxo(r)$ be bounded
in $r$, uniformly on compact subsets with respect to $\xo$.
We consider the matrix $P=(p_{ij})$ defined in polar coordinates 
$x=\xo+\rho e^{i\theta}$ by
\begin{equation}
\label{Pdef}
P(\xo+\rho e^{i\theta})=J(\theta)^*A(\xo+\rho e^{i\theta})J(\theta), 
\end{equation}
where $J$ is the rotation matrix defined in the Introduction.
Then, by \eqref{Jnabla} we have
\begin{align*}
\aij u_{x_i}u_{x_j}=&\bra A\nabla u,\nabla u\ket
=\bra AJ\onabla u,J\onabla u\ket=\bra P\onabla u,\onabla u\ket\\
=&p_{11}\urho^2+\left(p_{12}
+p_{21}\right)\urho\utr+p_{22}\left(\utr\right)^2.
\end{align*}
Note that $\det P=1$ and $P=P^*$. Therefore, denoting
$p_{11}=p$, $p_{12}=p_{21}=q$, $p_{22}=(1+q^2)/p$, we can write
\[
\aij u_{x_i}u_{x_j}=p\urho^2+2q\urho\utr+\frac{1+q^2}{p}\left(\utr\right)^2.
\]
In view of \eqref{elliptic} for every $k\in\R$ we have 
$\aij u_{x_i}u_{x_j}=((u-k)\aij u_{x_i})_{x_j}$
and therefore by the divergence theorem and an approximation argument (see \cite{PS})
\[
\gxo(r)=\int_\Sr(u-k)\aij u_{x_i}n_j\d\sigma
\]
where $\Sr$ is the circle of center $\xo$ and radius $r$, and 
$n$ is the outward normal to $\Sr$. Since $n=(\cos\theta,\sin\theta)=Je_1$
we have
\begin{align*}
\aij u_{x_i}n_j=\bra A\nabla u,n\ket
=\bra AJ\onabla u,Je_1\ket
=\bra P\onabla u,e_1\ket
=p\urho+q\utr
\end{align*}
and therefore by H\"older's inequality
\begin{align*}
\gxo(r)=&\int_\Sr(u-k)\left(p\urho+q\utr\right)\d\sigma
=\int_\Sr\sqrt{p}(u-k)\left(\sqrt{p}\urho
+\frac{q}{\sqrt{p}}\utr\right)\d\sigma
\\
\le&\left[\int_\Sr p(u-k)^2\d\sigma
\int_\Sr\left(\sqrt{p}\urho+\frac{q}{\sqrt{p}}\utr\right)^2\d\sigma\right]^{1/2}.
\end{align*}
Choosing $k=|\Sr|^{-1}\int_\Sr pu\d\sigma$, rescaling and using 
\eqref{wirtinger} with $a(\theta)=p(\xo+re^{i\theta})$,
$w(\theta)=u(\xo+re^{i\theta})-k$, we obtain the inequality
\[
\int_\Sr p(u-k)^2\d\sigma\le\left(|\Sr|^{-1}\int_\Sr p\d\sigma\right)^2
\int_\Sr\frac{1}{p}u_\theta^2\d\sigma.
\]
It follows that
\begin{align*}
\gxo(r)\le&\left(|\Sr|^{-1}\int_\Sr p\d\sigma\right)
\left[\int_\Sr\frac{1}{p}\utheta^2\d\sigma\right]^{1/2}
\left[\int_\Sr\left(\sqrt{p}\urho
+\frac{q}{\sqrt{p}}\utr\right)^2\d\sigma\right]^{1/2}.
\end{align*}
Equivalently, we may write
\begin{align*}
\gxo(r)\le&\left(|\Sr|^{-1}\int_\Sr p\d\sigma\right)r
\\
&\times\left[\int_\Sr\frac{1}{p}\left(\utr\right)^2\d\sigma\right]^{1/2}
\left[\int_\Sr\left(\sqrt{p}\urho
+\frac{q}{\sqrt{p}}\utr\right)^2\d\sigma\right]^{1/2}.
\end{align*}
By the inequality $\sqrt{ab}\le\frac{1}{2}(a+b)$ for every $a,b>0$ we obtain
for a.e.\ $0<r<\dist(\xo,\de\Omega)$ that
\begin{align}
\label{gineq}
\gxo(r)\le&\frac{1}{2}\left(|\Sr|^{-1}\int_\Sr p\d\sigma\right)r
\int_\Sr\left[\frac{1}{p}\left(\utr\right)^2
+\left(\sqrt{p}\urho
+\frac{q}{\sqrt{p}}\utr\right)^2\right]\d\sigma
\\
\nonumber
=&\frac{1}{2}\left(|\Sr|^{-1}\int_\Sr p\d\sigma\right)r
\int_\Sr\left[p\urho^2
+2q\urho\utr
+\left(\frac{1+q^2}{p}\right)\left(\utr\right)^2\right]\d\sigma
\\
\nonumber
=&\frac{1}{2}\left(|\Sr|^{-1}\int_\Sr p\d\sigma\right)r
\int_\Sr\bra P\onabla u,\onabla u\ket
\\
\nonumber
=&\frac{1}{2}\left(|\Sr|^{-1}\int_\Sr p\d\sigma\right)r
\int_\Sr\bra A\nabla u,\nabla u\ket\\
\nonumber
\le&\frac{1}{2}\es_{0<r<\dist(\xo,\de\Omega)}\left(|\Sr|^{-1}\int_\Sr p\d\sigma\right)
r\int_\Sr\bra A\nabla u,\nabla u\ket\\
\nonumber
\le&\frac{1}{2}\sup_{\xo\in\Omega}\es_{0<r<\dist(\xo,\de\Omega)}
\left(|\Sr|^{-1}\int_\Sr p\d\sigma\right)
r\int_\Sr\bra A\nabla u,\nabla u\ket.
\end{align}
Finally, we note that for all $\xi\in\R^2$, $|\xi|=1$, we have
$\xi=J(\theta)e_1$, $\theta=\arg\xi$, and therefore
\begin{align}
\label{p}
p(\xo+r\xi)=&\bra P(\xo+r\xi)e_1,e_1\ket=\bra J^*A(\xo+r\xi)J e_1,e_1\ket\\
\nonumber
=&\bra A(\xo+r\xi)\xi,\xi\ket
=a_{ij}(\xo+r\xi)\xi_i\xi_j.
\end{align}
Consequently, we have
\begin{align*}
\sup_{\xo\in\Omega}&\es_{0<r<\dist(\xo,\de\Omega)}
\left(|\Sr|^{-1}\int_\Sr p\d\sigma\right)\\
=&(2\pi)^{-1}\sup_{\xo\in\Omega}\es_{0<r<\dist(\xo,\de\Omega)}
\int_{|\xi|=1}a_{ij}(\xo+r\xi)\xi_i\xi_j\d\sigma\\
=&\alpha^{-1}.
\end{align*}
We note that $\gxo$ is differentiable almost everywhere and
that
\[
\gxo'(r)=\int_\Sr p\d\sigma\qquad\mathrm{a.e.}\ r.
\]
Therefore, \eqref{gineq} yields
\begin{equation}
\label{geq}
\gxo(r)\le\frac{r\gxo'(r)}{2\alpha}.
\end{equation}
Recalling that $G_{\xo}(r)=r^{-2\alpha}\gxo(r)$, we obtain
from \eqref{geq} that $(\log G_{\xo}(r))'\ge0$ a.e.
It follows that
\[
G_{\xo}(r)\le G_{\xo}(\dist(\xo,\de\Omega))
-\int_{r}^{\dist(\xo,\de\Omega)} G_{\xo}'(\rho)\,
\d\rho\le G_{\xo}(\dist(\xo,\de\Omega)).
\]
Hence, the desired H\"older estimate is established.
\end{proof}
\section*{Appendix}
We have postponed to this appendix a
\begin{proof}[Verification of Example~\ref{ex:sharp}]
In polar coordinates $x=\rho e^{i\theta}$, we define
$K(\theta)=\mathrm{diag}(k(\theta),1/k(\theta))$.
Then the matrix $A$ defined in \eqref{Adef} may be equivalently written in the form
\begin{align}
\label{Aexample}
A(x)=&J(\theta)K(\theta)J^*(\theta)\\
\nonumber
=&\left(
\begin{matrix}
&k(\theta)\cos^2\theta+\frac{1}{k(\theta)}\sin^2\theta
&&\left(k(\theta)-\frac{1}{k(\theta)}\right)\sin\theta\cos\theta\\
&\left(k(\theta)-\frac{1}{k(\theta)}\right)\sin\theta\cos\theta
&&k(\theta)\sin^2\theta+\frac{1}{k(\theta)}\cos^2\theta
\end{matrix}
\right),
\end{align} 
where
\[
J(\theta)=\left(
\begin{matrix}
&\cos\theta&&-\sin\theta\\
&\sin\theta&&\ \cos\theta
\end{matrix}\right).
\]
Hence, it is clear that $\det A\equiv1$. 
In order to check \eqref{exbalpha}, we note that if $\xo=0$ we have,
for all $0<r<1$:
\begin{align*}
a_{ij}(r\xi)\xi_i\xi_j
=\bra A(r\xi)\xi,\xi\ket=\bra J^*(\theta)A(r\xi)J(\theta)e_1,e_1\ket
=\bra K(\theta)e_1,e_1\ket=k(\theta),
\end{align*}
where $\theta=\arg\xi$.
Therefore,
\begin{equation}
\label{xo=0}
(2\pi)^{-1}\int_{|\xi|=1}a_{ij}(r\xi)\xi_i\xi_j\d\sigma_\xi
=(2\pi)^{-1}\int_0^{2\pi}k,
\end{equation}
for all $0<r<1$.
We assume that $k$ is smooth.
Then, for $\xo\neq0$, we have that $A$ is smooth near $\xo$.
It follows that
\begin{align*}
\inf_{0<\ro<\dist(\xo,\de B)}\es_{0<r<\ro}
\int_{|\xi|=1}a_{ij}(\xo+r\xi)\xi_i\xi_j\d\sigma_\xi
=\int_{|\xi|=1}a_{ij}(\xo)\xi_i\xi_j\d\sigma_\xi.
\end{align*}
In view of \eqref{Aexample}, we compute
\begin{align*}
\int_{|\xi|=1}a_{ij}(\xo)\xi_i\xi_j\d\sigma_\xi
=&\int_{|\xi|=1}\bra A(\xo)\xi,\xi\ket\d\sigma_\xi\\
=&\int_{|\xi|=1}\bra K(\thetao)J^*(\thetao)\xi,J^*(\thetao)\xi\ket\d\sigma_\xi\\
=&\int_0^{2\pi}\left\{k(\thetao)\cos^2(\varphi)+\frac{1}{k(\thetao)}\sin^2(\varphi)\right\}\d\varphi\\
=&\pi\left\{k(\thetao)+\frac{1}{k(\thetao)}\right\},
\end{align*}
where $\thetao=\arg\xo$.
Therefore, for every $\xo\neq0$ we obtain that
\begin{align}
\label{xonot0}
(2\pi)^{-1}\inf_{0<\ro<\dist(\xo,\de B)}\es_{0<r<\ro}
&\int_{|\xi|=1}a_{ij}(\xo+r\xi)\xi_i\xi_j\d\sigma_\xi\\
\nonumber
=&\frac{1}{2}\left\{k(\arg\xo)+\frac{1}{k(\arg\xo)}\right\}.
\end{align}
In view of \eqref{xo=0} and \eqref{xonot0}, formula \eqref{balpha}
takes the form:
\[
\balpha^{-1}=\max\left\{(2\pi)^{-1}\int_0^{2\pi}k,
\sup_{\xo\neq0}\frac{1}{2}\left(k(\arg\xo)+\frac{1}{k(\arg\xo)}\right)\right\}.
\]
Choosing $k(\theta)$ such that $M\le k\le 3M/2$,
for $M\gg1$ we achieve 
\[
\balpha=2\pi\left(\int_0^{2\pi}k\right)^{-1}.
\]
On the other hand,
it is readily checked that $u$ given by \eqref{uex} is a weak solution to the equation
\begin{equation}
\label{polar}
\left(\rho k(\theta)\urho\right)_\rho+\left(\frac{1}{k(\theta)}\frac{\utheta}{\rho}\right)_\theta=0,
\end{equation}
i.e.,
\begin{align*}
\int_B\left(k(\theta)\urho v_\rho+\frac{1}{k(\theta)}\utr
\frac{v_\theta}{\rho}\right)\rho\d\rho\d\theta=0
\end{align*}
for all $v\in H^1(B)$ compactly supported in $B$.
Recalling \eqref{Jnabla} and the definition of $K$,
we have
\begin{align*}
0=&\int_B\bra K\onabla u,\onabla u\ket\rho\d\rho\d\theta
=\int_B\bra KJ^*\nabla u,J^*\nabla u\ket\d x\\
=&\int_B\bra JKJ^*\nabla u,\nabla u\ket\d x
=\int_B\bra A\nabla u,\nabla u\ket\d x.
\end{align*}
Therefore, in cartesian coordinates, equation \eqref{polar} takes the form \eqref{elliptic}
with $A=(\aij )$ given by \eqref{Aexample}.
\end{proof}
\section*{Acknowledgements}
I should like to thank Professor Carlo Sbordone for
suggesting the problem, as well as for many interesting
and stimulating discussions.
I am also grateful to Professor Nicola Fusco for several useful comments on the manuscript.


\begin{thebibliography}{99}
\bibitem{DG}
E.~De Giorgi,
Sulla differenziabilit\`a e l'analiticit\`a delle
estremali degli integrali multipli, 
Mem.\ Acc.\ Sc.\ Torino \textbf{3} (1957), 25--42.
\bibitem{Gi}    
M.~Giaquinta,    
{\em Multiple Integrals in the Calculus of Variations
and Nonlinear Elliptic Systems},    
Annals of Math.\ Studies 105, Princeton Univ.\ Press, Princeton, 1983.    
\bibitem{GT}    
D.~Gilbarg and N.~Trudinger,    
{\em Elliptic Partial Differential Equations of Second Order},    
Classics in Mathematics, Springer-Verlag, Berlin/Heidelberg, 2001. 
\bibitem{IM}
T.~Iwaniec and G.~Martin,
\textit{Geometric Function Theory and Non-linear Analysis,}
Clarendon Press, Oxford, 2001.   
\bibitem{IS}
T.~Iwaniec and C.~Sbordone,
Quasiharmonic fields, 
Ann.\ I.H.P.\ \textbf{18} No.~5 (2001), 519--572.
\bibitem{Na}
J.~Nash,
Continuity of solutions of parabolic and elliptic equations, 
Amer.\ J.\  Math.\ \textbf{80} No.~5 (1958), 931--954.
\bibitem{PS} 
L.C.~Piccinini and S.~Spagnolo, 
On the H\"older continuity of solutions of second order elliptic equations
in two variables, Ann.\ Scuola Norm.\ Sup.\ Pisa \textbf{26} No.~2 (1972), 391--402.
\end{thebibliography}
\end{document}